\documentclass[11pt, oneside]{amsart}
\usepackage{amstext}
\usepackage{amsmath, amssymb, amsthm}
\usepackage{amssymb}
\usepackage{mathrsfs}
\usepackage{enumerate}
\usepackage{xcolor}
\usepackage{comment}
\usepackage{hyperref}
\hypersetup{
    colorlinks=true,    % false: boxed links; true: colored links
    linkcolor=blue,     % Color of internal links (sections, TOC, etc.)
    citecolor=blue,    % Color of citation links (bibliography)
    filecolor=blue,  % Color of local file links
    urlcolor=blue       % Color of web URLs
}
\usepackage{graphicx}
\usepackage{mathtools}
\usepackage{accents}
\mathtoolsset{showonlyrefs}
\usepackage[a4paper, total={6.5in, 8.6in}]{geometry}
\usepackage{setspace}

\newtheorem{theorem}{Theorem}[section]
\newtheorem{lemma}[theorem]{Lemma}
\newtheorem{proposition}[theorem]{Proposition}
\newtheorem{corollary}[theorem]{Corollary}

\newtheorem{question}[theorem]{Question}

\newtheorem{example}[theorem]{Example}

\DeclareMathOperator{\Ric}{Ric}
\DeclareMathOperator{\dvol}{dvol}

\DeclareMathOperator{\dist}{dist}

\title{The Cheng--Yau gradient estimate for Jacobi fields}
\author{Stephen Lynch}
\address{Department of Mathematics, King's College London, Strand, London, WC2R 2LS, UK}
\email{stephen.lynch@kcl.ac.uk}
\date{}

\setdisplayskipstretch{1}

\begin{document}

\begin{abstract}
We prove a pointwise gradient estimate for positive Jacobi fields on minimal hypersurfaces of dimension two and three, inspired by the Cheng--Yau estimate for harmonic functions. The estimate also controls the second fundamental form, and so provides yet another proof of the stable Bernstein theorem in dimensions two and three. 
\end{abstract}

\maketitle

\noindent\textbf{Use of AI.} This project grew out of a question posed by the author to an LLM: is there an analogue of the Cheng--Yau gradient estimate for positive Jacobi fields on a minimal surface? The model quickly found such an estimate, gave a nice proof, and observed that it also implies a curvature estimate. After more prompting the model found a generalisation to dimension three, which was substantially weaker in that the oscillation of the Jacobi field appeared on the right-hand side, and gave a fairly incomprehensible proof. With some work, and further input from the LLM, we were able to remove the dependence on the oscillation using the conformally changed metric of Catino--Mastrolia--Roncoroni~\cite{CatinoMastroliaRoncoroni}, and to find a transparent proof using an idea of Colding~\cite{Colding}. The LLM was also used for mathematical discussion, checking, and editorial feedback. The author takes full responsibility for the final manuscript.

\section{Introduction}

One of the basic principles of geometric analysis is that, in the presence of lower curvature bounds, positive solutions to elliptic and parabolic PDE should not oscillate too rapidly. While there are numerous examples of this phenomenon, surely one of the most beautiful is the interior gradient estimate of Cheng--Yau~\cite{ChengYau} (see also \cite{Yau}): given a positive harmonic function $u$, defined on a complete Riemannian $n$-manifold $M$ with Ricci curvature at least $-(n-1)K$ for some $K \geq 0$, we have 
\[
|\nabla \log u(p)| \leq (n-1)\sqrt{K} + \frac{C}{\dist(p,\partial M)}
\]
at every point $p \in M \setminus \partial M$. The constant $C$ depends only on $n$. The proof of the Cheng--Yau estimate is based on the maximum principle and a sharp lower bound for $\Delta|\nabla \log u|^2$ obtained from the Bochner formula. Its influence in geometry has been far-reaching and profound. 

In this paper we derive an analogue of the Cheng--Yau estimate on minimal hypersurfaces in $\mathbb{R}^3$ and $\mathbb{R}^4$. Minimal hypersurfaces clearly do not have Ricci curvature bounded below a priori, but by the Gauss equation and minimality we have the simple identity
\[
\Ric = -A^2,
\]
where $A$ is the (scalar-valued) second fundamental form. This suggests that instead of harmonic functions we look at solutions of an equation which couples to curvature. A natural candidate is the Jacobi equation,
\[
\Delta u + |A|^2 u =0.
\]
Solutions to this equation are called Jacobi fields. They play an important role in the study of minimal hypersurfaces, because they form the kernel of the second variation of area. 

Our main result is the following interior gradient estimate for positive Jacobi fields. It bounds both $|\nabla\log u|^2$ and $|A|^2$, with a constant independent of the size or oscillation of $u$, and independent of the volume of the minimal hypersurface. To prove the estimate we combine the Bochner formula for $\Delta |\nabla \log u|^2$ with Simons' identity~\cite{Simons} and apply the maximum principle.

\begin{theorem}\label{main}
Assume $n \in \{2,3\}$. Let $\Sigma \to \mathbb{R}^{n+1}$ be a complete minimally immersed hypersurface, equipped with a Jacobi field $u$ which is positive on $\Sigma\setminus\partial\Sigma$. We then have 
\[
|\nabla \log u(p)|^2 + |A(p)|^2 \leq C\dist(p,\partial\Sigma)^{-2}
\]
at every point $p \in \Sigma\setminus \partial\Sigma$, where $C$ is a universal constant and $\dist$ is intrinsic. 
\end{theorem}

There is a clear connection between this estimate and the curvature estimates for minimal graphs due to Caffarelli--Nirenberg--Spruck \cite{CaffarelliNirenbergSpruck} and Ecker--Huisken \cite{EckerHuisken}. These are proved by bounding a suitable combination of $|A|^2$ and the vertical component of the unit normal $\langle \nu, e\rangle$; the latter quantity is, in the graphical setting, a positive Jacobi field.

Integrating the gradient estimate along geodesics gives a Harnack inequality.

\begin{corollary}\label{Harnack}
Let $\Sigma$ and $u$ be as in Theorem~\ref{main}, and let $C$ denote the constant in that theorem. For every $p\in\Sigma\setminus\partial\Sigma$ and $0<r<R\leq\dist(p,\partial\Sigma)$, we have
\[
\frac{\sup_{B(p,r)}u}{\inf_{B(p,r)}u}
\leq
\left(\frac{R}{R-r}\right)^{2\sqrt{C}},
\]
where $B(p,r)$ is the intrinsic ball of radius $r$ about $p$.
\end{corollary}

Harnack inequalities on minimal hypersurfaces go back to Bombieri--Giusti~\cite{BombieriGiusti}.  If curvature is already controlled, the estimate in Corollary~\ref{Harnack} can be derived from standard elliptic theory. The point here is that we obtain a curvature estimate and the Harnack inequality together, from a single application of the maximum principle. 

%Gradient estimates for minimal graphs were proved in \cite{BombieriDeGiorgiMiranda}, \cite{Trudinger} and \cite{Korevaar}.

\subsection{Bernstein theorems} For two-sided minimal hypersurfaces, the existence of a positive Jacobi field is equivalent to stability~\cite{FischerColbrieSchoen}. Thus Theorem~\ref{main} recovers some important results concerning stable minimal hypersurfaces.

Dropping the gradient term in Theorem~\ref{main} gives a new proof of the following curvature estimate, which is due to Schoen~\cite{Schoen} in dimension two and Chodosh--Li~\cite{ChodoshLi} in dimension three. 

\begin{corollary}\label{curvature est}
Assume $n \in \{2,3\}$. Let $\Sigma\to\mathbb{R}^{n+1}$ be a complete minimally immersed hypersurface which is two-sided and stable. We then have
\[
|A(p)|^2\leq C\dist(p,\partial\Sigma)^{-2}
\]
at every interior point, where $C$ is universal.
\end{corollary}

Applying this estimate on larger and larger intrinsic balls gives the stable Bernstein theorem. 

\begin{corollary}\label{Bernstein}
For $n\in\{2,3\}$, the only complete, connected, two-sided, stable minimal immersions without boundary in $\mathbb{R}^{n+1}$ are hyperplanes.
\end{corollary}

In dimension two, Bernstein~\cite{Bernstein} proved that entire minimal graphs are planes in 1915. This was generalised to complete, stable, two-sided minimal immersions by do Carmo--Peng~\cite{doCarmoPeng}, Fischer-Colbrie--Schoen~\cite{FischerColbrieSchoen} and Pogorelov~\cite{Pogorelov} around 1980.

In dimensions $3\leq n\leq6$, Bernstein rigidity and curvature estimates were established under natural volume bounds by Schoen--Simon--Yau~\cite{SchoenSimonYau}, Schoen--Simon~\cite{SchoenSimon} and Bellettini~\cite{Bellettini}, whereas counterexamples exist for $n\geq7$~\cite{BombieriDeGiorgiGiusti, HardtSimon}.

Without a volume growth assumption, Chodosh--Li~\cite{ChodoshLi} only recently settled the stable Bernstein problem in dimension $n=3$. They used a monotonicity formula for the level sets of the Green's function, and later gave a $\mu$-bubble proof~\cite{ChodoshLiAnisotropic}. Another approach was found by Catino--Mastrolia--Roncoroni~\cite{CatinoMastroliaRoncoroni}, using weighted volume comparison for the conformal metric $u^{4/3}g$ (this metric also appears in our arguments). More recently, Cabr\'{e}--Catino--Mari--Mastrolia--Roncoroni~\cite{CabreCatinoMariMastroliaRoncoroni} found a further alternative proof using a gradient estimate for the Green's function. These approaches obtain flatness from integral curvature estimates, by establishing volume bounds or using cutoffs built from the Green's function.

The $\mu$-bubble approach was extended to dimension $n=4$ by Chodosh--Li--Minter--Stryker~\cite{ChodoshLiMinterStryker}, and to $n =5$ by Mazet~\cite{Mazet} (see also \cite{AntonelliXu}), establishing the stable Bernstein theorem in these dimensions. The article \cite{HongLiWang} deals with the final remaining dimension, $n=6$.

\subsection{The proof} Let us outline the key steps in the proof of Theorem~\ref{main}. We focus on the case $n=3$, since it is more involved and implies the two-dimensional result. 

First we recall a remarkable fact that has often been used in the study of stable minimal hypersurfaces: given a positive Jacobi field $u$, one can weight the induced metric $g$ by a suitable power of $u$ so that it acquires nonnegative scalar curvature.\footnote{This goes back at least to the work of Fischer-Colbrie~\cite{FischerColbrie} and, as a broader principle, the conformal descent scheme of Schoen--Yau~\cite{SchoenYauStructure}.} For $n =3$ the calculation goes as follows. Setting $h = u^{2k} g$, the standard conformal scalar curvature formula and $R_g = -|A|^2$ (which follows from minimality and the Gauss equation) yield
\iffalse
\[
u^{2k} R_h = R_g  + 4k|A|^2 + 2k(2-k)|\nabla \log u|^2.
\]
By minimality and the Gauss equation we have $R_g = -|A|^2$, so this becomes
\fi
\begin{equation}\label{conformal R intro}
u^{2k} R_h = (4k-1)|A|^2 + 2k(2-k)|\nabla \log u|^2.
\end{equation}
Therefore, $R_h$ is nonnegative for $\frac{1}{4} \leq k \leq 2$.

The right-hand side of \eqref{conformal R intro} contains precisely the two quantities we must estimate to establish Theorem~\ref{main}. This suggests that we derive an upper bound for the scalar curvature $R_h$. We choose $k=\frac{2}{3}$, following Catino--Mastrolia--Roncoroni~\cite{CatinoMastroliaRoncoroni}, who showed that for this choice $h$ has nonnegative $5$-dimensional Bakry--\'{E}mery Ricci curvature. The formula~\eqref{conformal R intro} then becomes
\begin{equation}\label{conformal R intro 2/3}
R_h = u^{-4/3}\bigg(\frac53|A|^2+\frac{16}{9}|\nabla\log u|^2\bigg).
\end{equation}

The most important ingredient in the proof of Theorem~\ref{main} is the inequality
\begin{equation}\label{Delta R intro}
\Delta_h R_h \geq \frac{2}{75} R_h^2,
\end{equation}
which we derive from the Bochner formula and Simons' identity; more is said about the calculation just below. It is well known that functions which satisfy an inequality of this form enjoy strong local estimates. Assuming $\Sigma$ is complete, we use \eqref{Delta R intro} and the maximum principle to bound $R_h$ at the centre of any $h$-ball in terms of its radius. This requires controlled cutoffs, which are available because of the Laplacian comparison theorem for metrics with nonnegative finite-dimensional Bakry--\'{E}mery Ricci curvature~\cite{Qian}. The remaining step is to convert our bound for $R_h$ in $h$-balls into the estimate claimed in Theorem~\ref{main}. To do this we prove a lemma which compares $h$-distances with $g$-distances and use \eqref{conformal R intro 2/3}. The distance comparison absorbs the weight $u^{-4/3}$ in \eqref{conformal R intro 2/3}, so that we obtain an estimate for the quantity $\frac{5}{3}|A|^2 + \frac{16}{9}|\nabla \log u|^2$ which immediately implies Theorem~\ref{main}. 

\subsubsection{Combining Bochner and Simons.} Let us describe the basic mechanism that leads to \eqref{Delta R intro}. It is instructive to ask what happens if we remain in the induced metric and try to apply the maximum principle to a linear combination of $|A|^2$ and $|\nabla \log u|^2$. Writing $f=\log u$, the Jacobi equation gives
\begin{equation}\label{Delt f intro}
\Delta f=-|A|^2-|\nabla f|^2.
\end{equation}
Combining this with the Bochner formula and $\Ric = - A^2$, we obtain
\[
\Delta |\nabla f|^2 = 2|\nabla^2 f|^2 - 2\langle\nabla f, \nabla |A|^2 + \nabla |\nabla f|^2\rangle -  2|A(\nabla f, \cdot)|^2.
\]
At the same time, Simons' identity asserts that 
\[
\Delta |A|^2 = 2|\nabla A|^2 - 2|A|^4. 
\]
This coupled system of equations has the following useful features:
\begin{itemize}
    \item By \eqref{Delt f intro}, the Hessian term in the Bochner formula
    contains favourable terms which are quartic in $|\nabla f|$ and $|A|$.\footnote{This is very similar to the structure exploited by Cheng--Yau~\cite{ChengYau} and Yau~\cite{Yau} to prove their estimates for positive harmonic functions.} These are already visible in the trace decomposition
\[
|\nabla^2 f|^2
= \left|\nabla^2 f-\frac{1}{3}\Delta f\,g\right|^2
+ \frac{1}{3}(|\nabla f|^2+|A|^2)^2.
\]
    They help counteract both the negative curvature term in Simons'
    identity and the Ricci term $-2|A(\nabla f,\cdot)|^2$.
    \item Simultaneously, the favourable gradient term in Simons' identity can be used to combat the cross-term $-2\langle\nabla f, \nabla |A|^2\rangle$ in the Bochner formula, using the sharp Kato-type inequality 
\[
|\nabla A|^2\geq\frac5{12}\frac{|\nabla|A|^2|^2}{|A|^2}
\]
    of Schoen--Simon--Yau \cite{SchoenSimonYau}.
\end{itemize}

For a general linear combination $Q=c|A|^2 + |\nabla f|^2$ with $c > 0$, by combining Bochner and Simons and exploiting the various favourable terms that arise in the manner outlined above, one can derive the inequality
\begin{equation}\label{Delta Q general}
\Delta Q\geq
\frac34\frac{|\nabla Q|^2}{Q}
-\langle\nabla f,\nabla Q\rangle
+|\nabla f|^4
+\left(\frac{20}{3}-3c-\frac{10}{3c}\right)|\nabla f|^2|A|^2
+\left(\frac23-2c\right)|A|^4.
\end{equation}
Unfortunately there is no choice of $c$ which renders both of the coefficients on the right nonnegative. 

Passing to the conformally changed metric $h = u^{4/3} g$, and treating the conformal scalar curvature $R_h$, improves the situation. By \eqref{conformal R intro 2/3}, with $c = \frac{15}{16}$ we have $R_h = \frac{16}{9}u^{-4/3} Q$, and the standard formula for the Laplacian under a conformal change gives
\[
u^{8/3}\Delta_h(u^{-4/3}Q) = \Delta Q - 2\langle\nabla f,\nabla Q\rangle + \bigg(\frac{20}{9}|\nabla f|^2 + \frac{4}{3}|A|^2\bigg)Q.
\]
Notice that we have gained additional favourable terms on the right---these arise when the Laplacian falls on the weight $u^{-4/3}$, because of the Jacobi equation. It is these additional terms, combined with the Bochner--Simons calculation described above, that lead to the key inequality \eqref{Delta R intro}. 

Curiously, we do not need to make any conformal change to prove Theorem~\ref{main} in dimension two. There the gradient term in Simons' identity is more powerful, and the negative Ricci curvature of $\Sigma$ is smaller relative to $|A|^2$, with the consequence that our maximum principle argument closes in the original metric $g$. This simpler proof is included in Section~\ref{sec 2D proof}.

\subsection{Outline} In Section~\ref{sec BS} we carry out the Bochner--Simons calculation which gives \eqref{Delta R intro}. We also the model case where $\Sigma$ is flat and $u$ is a positive harmonic function; this is intended to make the geometry underlying \eqref{Delta R intro} more apparent. Theorem~\ref{main} is then proved in Section~\ref{sec main proof}. In Section~\ref{sec 2D proof} we give an alternative proof which only works when $n = 2$, but is a lot simpler. In Appendix~\ref{sharp BS} we prove a a refined version of \eqref{Delta R intro}, which is not needed to prove our main theorem, but may be of interest in its own right. To make this article self-contained we have included Appendix~\ref{Laplace comparison}, which recalls the weighted Laplacian comparison estimate from \cite{Qian} used in Section~\ref{sec main proof}, and Appendix~\ref{conformal formulae}, which records some standard conformal transformation formulae.  

\subsection{Further questions} The following questions are left open by this work.

\begin{question}[Higher dimensions] For $4 \leq n \leq 6$, can one prove an estimate as in Theorem~\ref{main} using the maximum principle?
\end{question}

One difficulty seems to be that the Gauss equation and minimality provide less control on the Ricci curvature in higher dimensions. When $n \geq 4$, for certain powers $k$ and $h = u^{2k}g$ the scalar curvature $R_h$ does satisfy a nice elliptic equation, but we have not found an inequality like \eqref{Delta R intro}. One could experiment with combining $|\nabla \log u|^2$ and $|A|^2$ in different ways.

\begin{question}[Sharp constant] What is the sharp constant in Theorem~\ref{main}, and which pairs $\Sigma$ and $u$ realise it?
\end{question}

Tracking the constants in our proof shows that we can take $C = 10^{47}$, but this is surely far from optimal. 

\begin{question}[Parabolic analogues]
What is the natural analogue of Theorem~\ref{main} for positive solutions to the equation
\[
\partial_t u = \Delta u + |A|^2 u
\]
along a solution to the mean curvature flow?
\end{question}

Here we are motivated by the Li--Yau inequality~\cite{LiYau} and Hamilton's estimate~\cite{Hamilton}.

\section{The Bochner--Simons calculation in dimension three}\label{sec BS}

Consider a minimal hypersurface $\Sigma$ immersed in $\mathbb{R}^4$ which carries a positive Jacobi field $u$. Let us set
\[
f=\log u,\qquad q=|\nabla f|^2,\qquad s=|A|^2,
\qquad h=u^{4/3}g.
\]
All unadorned derivatives and norms are taken with respect to $g$. The Jacobi equation and the Gauss equation give
\[
\Delta f=-q-s,\qquad \Ric=-A^2,\qquad R=-s.
\]
Combining this with the standard conformal scalar curvature formula shows that
\begin{equation}\label{Jacobi scalar curvature}
R_h=u^{-4/3}\left(R-\frac{8}{3}\Delta f-\frac{8}{9}q\right)
=u^{-4/3}\left(\frac{16}{9}q+\frac{5}{3}s\right).
\end{equation}

Using the Bochner formula and Simons' identity, we will demonstrate that $R_h$ is a subsolution to an elliptic equation with quadratic nonlinearity. As a warm-up, it is useful to consider the case $\Sigma = \mathbb{R}^3$. 

\begin{example}\label{harmonic example}
On flat space, the Jacobi equation reduces to $\Delta u = 0$. The model for our calculations in this section is the case where $u$ is the fundamental solution $u_* = r^{-1}$ defined on $\mathbb{R}^3 \setminus \{0\}$. The corresponding metric $h_* = u^{4/3}_* g_{\mathbb{R}^3}$ is a cone; indeed, writing the Euclidean metric in polar coordinates and setting $\rho=3r^{1/3}$, we obtain
\[
h_*=r^{-4/3}(dr^2+r^2g_{\mathbb{S}^2})
=d\rho^2+\frac{\rho^2}{9}g_{\mathbb{S}^2}.
\]
Moreover, in this case the scalar curvature satisfies
\begin{equation}\label{cone scalar curvature}
R_{h_*}=\frac{16}{\rho^2},
\qquad
\Delta_{h_*} R_{h_*}=\frac{1}{8}R_{h_*}^2.
\end{equation}

Up to a nonnegative defect, the second identity in \eqref{cone scalar curvature} persists for an arbitrary positive harmonic function. The defect is measured by the traceless tensor
\begin{equation}\label{T harmonic}
T=-\frac{u^2}{2}\big(\nabla^2(u^{-2})\big)^o
=\nabla^2f-2df\otimes df+qg,
\end{equation}
which played a prominent role in Colding's work~\cite{Colding}. This tensor vanishes on the fundamental solution, since the Hessian of $u_*^{-2} = r^2$ is pure trace. Conversely, one can show that $T$ vanishes identically if and only if $u$ is the fundamental solution (up to scaling) or a constant.

Inserting $\Delta f = -q$ into the Bochner formula, and expressing $\nabla^2 f$ in terms of $T$, one obtains
\begin{equation}\label{Bochner harmonic T}
\Delta q-2\langle\nabla f,\nabla q\rangle=2|T|^2-2q^2.
\end{equation}
On the other hand, standard conformal transformation formulae yield
\begin{align*}
\Delta_h R_h
&=\frac{16}{9}u^{-8/3}
\left(\Delta q-2\langle\nabla f,\nabla q\rangle+\frac{20}{9}q^2\right).
\end{align*}
Together these identities imply 
\begin{equation}\label{Delta R harmonic}
\Delta_h R_h - \frac18R_h^2=\frac{32}{9}u^{-8/3}|T|^2.
\end{equation}

In summary, for any positive harmonic function defined on $\mathbb{R}^3\setminus\{0\}$, the difference $\Delta_h R_h - \frac{1}{8} R_h^2$ is nonnegative, and it vanishes identically precisely when $u = u_*$ (up to scaling) or $u$ is constant.
\end{example}

Much of the structure described in Example~\ref{harmonic example} finds analogy in our more general setting, where $\Sigma$ is not necessarily flat but is minimal and stable. We still define $T$ to be the traceless part of $-\frac{u^2}{2}\big(\nabla^2(u^{-2})\big)$, but now since $\Delta u = -su$ this means setting
\begin{equation}\label{T general}
T =\nabla^2f-2df\otimes df+\left(q+\frac{s}{3}\right)g.
\end{equation}
Let us mention again that our introduction of this tensor is inspired by~\cite{Colding}.%\footnote{We originally used the tensor $T$ purely as a tool to extract information from the Hessian term in the Bochner formula (see the proof of Proposition~\ref{Jacobi scalar inequality}). In asking why this works well, we were led to the discussion in Example~\ref{harmonic example} and the geometric interpretation of Theorem~\ref{main} (put forward in the introduction) as a scalar curvature estimate for the metric $h$. This interpretation greatly influenced the version of the proof finally recorded here.}

The inequality derived in the following proposition is sufficient to prove Theorem~\ref{main}. For the sake of completeness, in Appendix~\ref{sharp BS} we have included a refinement which throws nothing away and recovers \eqref{Delta R harmonic} on flat space. 

\begin{proposition}\label{Jacobi scalar inequality}
The scalar curvature of $h=u^{4/3}g$ satisfies
\begin{equation}\label{Delta scalar curvature}
\Delta_h R_h\geq\frac{2}{75}R_h^2.
\end{equation}
\end{proposition}
\begin{proof}
Combining the Bochner formula with $\Delta f=-q-s$ and $\Ric=-A^2$, we obtain
\[
\Delta q=2|\nabla^2f|^2
-2\langle\nabla f,\nabla q+\nabla s\rangle
-2|A(\nabla f,\cdot)|^2.
\]
Rewriting the Hessian term using
\[
|T|^2=|\nabla^2f|^2
-2\langle\nabla f,\nabla q\rangle
+q^2-2qs-\frac13s^2
\]
then gives
\begin{align}\label{Bochner Jacobi T}
\Delta q &= 2|T|^2+2\langle\nabla f,\nabla q - \nabla s\rangle
-2|A(\nabla f,\cdot)|^2\notag\\
&\qquad -2q^2+4qs+\frac23s^2.
\end{align}
We combine this with Simons' identity
\[
\Delta s=2|\nabla A|^2-2s^2
\]
to conclude that, for any constant $c$, 
\begin{align}\label{Delta q + cs}
\Delta (q+cs) &= 2|T|^2 + 2c|\nabla A|^2 +2\langle\nabla f,\nabla q - \nabla s\rangle \notag\\
&\qquad -2|A(\nabla f,\cdot)|^2 -2q^2+4qs+\left(\frac{2}{3} - 2c\right)s^2. 
\end{align}

The formula for the Laplacian under a conformal change gives
\[
u^{4/3}\Delta_h = \Delta + \frac{2}{3}\langle\nabla f, \nabla \cdot\rangle.
\]
After a short calculation using $\Delta u = -su$ one finds that, for an arbitrary function $Q$,
\begin{equation}\label{weighted conformal Laplacian}
u^{8/3}\Delta_h(u^{-4/3}Q)
=\Delta Q-2\langle\nabla f,\nabla Q\rangle
+\left(\frac{20}{9}q+\frac43s\right)Q. 
\end{equation}
We set $Q = \frac{16}{9}q+\frac{5}{3}s$, so that $u^{4/3} R_h = Q$, and combine \eqref{weighted conformal Laplacian} and \eqref{Delta q + cs} with $c = \frac{15}{16}$. This gives
\begin{align}\label{Bochner Simons scalar}
u^{8/3}\Delta_h R_h
={}&\frac{32}{9}|T|^2+\frac{10}{3}|\nabla A|^2
-\frac{62}{9}\langle\nabla f,\nabla s\rangle\notag\\
&-\frac{32}{9}|A(\nabla f,\cdot)|^2
+\frac{32}{81}q^2+\frac{356}{27}qs+\frac{2}{27}s^2.
\end{align}

At a point where $s=0$, we have $A=0$ and $\nabla s=0$, and in this case \eqref{Delta scalar curvature} follows directly from \eqref{Bochner Simons scalar}. We are free therefore to assume that $s>0$. We use the Kato-type inequality of Schoen--Simon--Yau~\cite{SchoenSimonYau},
\begin{equation}\label{Kato}
|\nabla A|^2\geq\frac5{12}\frac{|\nabla s|^2}{s},
\end{equation}
to absorb the cross-term in \eqref{Bochner Simons scalar} involving $\nabla s$: completing the square shows that
\begin{align*}
\frac{10}{3}|\nabla A|^2
-\frac{62}{9}\langle\nabla f,\nabla s\rangle
&\geq\frac{25}{18s}
\left|\nabla s-\frac{62}{25}s\nabla f\right|^2
-\frac{1922}{225}sq.
\end{align*}
Since $A$ is traceless and we are in dimension three, we also have
\[
|A(\nabla f,\cdot)|^2\leq\frac23sq.
\]
Inserting these bounds into \eqref{Bochner Simons scalar} gives
\begin{align*}
u^{8/3}\Delta_h R_h
&\geq\frac{32}{9}|T|^2
+\frac{32}{81}q^2+\frac{1534}{675}qs+\frac{2}{27}s^2
\end{align*}
and hence, after discarding $|T|^2$ and part of the reaction terms,
\[
u^{8/3}\Delta_h R_h \geq \frac{2}{75} \left(\frac{16}{9}q+\frac{5}{3}s\right)^2. 
\]
This proves \eqref{Delta scalar curvature}.
\end{proof}

\section{Application of the maximum principle}\label{sec main proof}

We continue in the setting of Section~\ref{sec BS}, but now additionally assume that $\Sigma$ is complete. Our next theorem exploits the favourable quadratic reaction term in \eqref{Delta scalar curvature} to bound $R_h$, using the maximum principle and cutoffs constructed from the $h$-distance function. To control the cutoffs, we will need the following Ricci curvature estimate from \cite{CatinoMastroliaRoncoroni}: writing $\psi = \frac{2}{3}f$, we have
\begin{equation}\label{BE}
\Ric_h+ \nabla^2_h\psi-\frac{1}{2}d\psi\otimes d\psi \geq 0.
\end{equation}
To verify this, one uses the standard formula for the Ricci curvature under a conformal change to express
\begin{equation}\label{conformal Ricci}
\Ric_h
=
\Ric-\frac{2}{3}\nabla^2_hf
-\frac{4}{9}df\otimes df
-\frac{2}{3}\Delta f\,g,
\end{equation}
and then inserts
\[
\Ric=-A^2\geq-\frac{2}{3}sg,
\qquad
df\otimes df\leq qg,
\qquad
\Delta f=-s-q,
\]
where the first inequality holds because $n = 3$ and $A$ is traceless. 
\iffalse
The Ricci and Hessian formulae for a conformal metric $e^{2\phi}g$ in dimension three are
\begin{align}\label{Ricci h}
\Ric_h
&=
\Ric-\nabla^2\phi+d\phi\otimes d\phi
-(\Delta\phi+|\nabla\phi|^2)g,\\
\nabla^2_hv
&=
\nabla^2v-dv\otimes d\phi-d\phi\otimes dv
+\langle\nabla v,\nabla\phi\rangle g. \notag
\end{align}
Substituting the second formula into the first, and setting $\phi=\frac{k}{2}f$, we obtain
\[
\Ric_h
=
\Ric-\frac{k}{2}\nabla^2_hf
-\frac{k^2}{4}df\otimes df
-\frac{k}{2}\Delta f\,g.
\]
We now combine this with
\[
\Ric=-A^2\geq-\frac{2}{3}sg,
\qquad
df\otimes df\leq qg,
\qquad
\Delta f=-s-q,
\]
and so arrive at the lower bound
\[
\Ric_h
\geq
-\frac{k}{2}\nabla^2_hf
-\left(\frac{k^2}{4}-\frac{k}{2}\right)df\otimes df
+\left(\frac{k}{2}-\frac{2}{3}\right)sg.
\]
This shows why our choice of $k=\frac{4}{3}$ is convenient---it cancels the final term on the right. Introducing the notation
\[
\psi=\frac{2}{3}f,
\]
we are then left with
\[
\Ric_h+\nabla^2_h\psi-\frac{1}{2}d\psi\otimes d\psi\geq0.
\]
\fi
The tensor on the left of \eqref{BE} is the $5$-dimensional Bakry--\'{E}mery Ricci tensor for the weighted manifold $(\Sigma,h,e^{-\psi}\dvol_h)$.

\begin{samepage}
\begin{theorem}\label{conf max princ}
If $\overline{B_h(p,R)}$ is compactly contained in $\Sigma \setminus \partial\Sigma$, then
\[
R_h(p)\leq KR^{-2},
\]
where $K$ is a universal positive constant.
\end{theorem}
\end{samepage}
\begin{proof}
We introduce a nonnegative cutoff function $\varphi$ which vanishes on the boundary of $\overline{B_h(p,R)}$. If $R_h$ vanishes identically then there is nothing to prove, and otherwise there is a point $\bar p$ where $\varphi R_h$ attains a positive interior maximum. At $\bar p$ we have
\[
\nabla^h(\varphi R_h)=0,
\qquad
\Delta_h(\varphi R_h)\leq0.
\]
\iffalse
Combining these with the formulae
\[
\nabla^h(\varphi W)
=
W\nabla^h\varphi+\varphi\nabla^h W
\]
and
\[
\Delta_h(\varphi W)
=
W\Delta_h\varphi
+\varphi\Delta_h W
-2\varphi^{-1}|\nabla^h\varphi|_h^2W
+2\varphi^{-1}
\langle\nabla^h\varphi,
\nabla^h(\varphi W)\rangle_h,
\]
we deduce that the inequality
\begin{equation}\label{max phi W}
\left(
\Delta_h\varphi
-2\varphi^{-1}|\nabla^h\varphi|_h^2
\right)W
+\varphi\Delta_h W
\leq0
\end{equation}
holds at $\bar p$. Discarding the nonnegative gradient term in \eqref{hat Delta W square}, we obtain
\[
\Delta_h W\geq\frac{2}{9}W^2.
\]
\fi
After a straightforward calculation using \eqref{Delta scalar curvature}, these conditions lead to
\begin{equation}\label{max phi W simple}
\Delta_h\varphi
-2\varphi^{-1}|\nabla^h\varphi|_h^2
+\frac{2}{75}\varphi R_h
\leq0.
\end{equation}

Let us set $\zeta(x)=R^2-r(x)^2$, where $r$ is the $h$-distance function from $p$, and define $\varphi(x)=\zeta(x)^2$. We then have
\[
\varphi^{-1}|\nabla^h\varphi|_h^2=16r^2, 
\qquad 
\Delta_h\varphi=-4\zeta(1+r\Delta_h r)+8r^2,
\]
and hence, from \eqref{max phi W simple},
\begin{equation}\label{max phi W II}
-4\zeta(1+r\Delta_h r)
-24r^2
+\frac{2}{75}\varphi R_h
\leq 0.
\end{equation}
At a cut point this calculation is made rigorous by approximation of $r$ through smooth upper barriers (see Appendix~\ref{Laplace comparison}). 

It remains to estimate $r\Delta_h r$ from above. Nonnegativity of the $5$-dimensional Bakry--\'{E}mery Ricci tensor (see \eqref{BE}) furnishes us with the comparison estimate\footnote{This estimate is essentially contained in \cite{Qian}, but we provide a proof in Appendix~\ref{Laplace comparison} for the sake of completeness.}
\[
\Delta_{h,\psi} r\leq\frac{4}{r},
\]
where $\Delta_{h,\psi}$ is the drift Laplacian $\Delta_h -\langle\nabla^h\psi, \nabla^h\cdot\rangle_h$. Inserting this estimate into \eqref{max phi W II} shows that, at the point $\bar p$,
\[
-20\zeta
-24r^2
-4r\zeta
\langle\nabla^h\psi,\nabla^h r\rangle_h
+\frac{2}{75}\varphi R_h
\leq0.
\]

We estimate the inner-product term using
\[
|\nabla^h\psi|_h^2
=
\frac{4}{9}e^{-\frac{4}{3}f}q
\leq
\frac{1}{4}R_h
\]
and Cauchy--Schwarz. This gives
\[
-20\zeta
-24r^2
-2r\sqrt{\varphi R_h}
+\frac{2}{75}\varphi R_h
\leq0
\]
at the point $\bar p$. Inserting $20\zeta+24r^2=20R^2+4r^2\leq24R^2$, and using Young's inequality to absorb the square-root term, we finally arrive at 
\[
\varphi(\bar p) R_h(\bar p) \leq 7500R^2.
\]
The claim now follows, since maximality of $\bar p$ implies
\[
R^4 R_h(p) = \varphi(p)R_h(p) \leq \varphi(\bar p) R_h(\bar p) \leq 7500R^2. \qedhere
\]
\end{proof}

Let us define
\[
Q=\frac{16}{9}q+\frac53s,
\]
so that $R_h=u^{-4/3}Q$. Theorem~\ref{conf max princ} bounds $R_h$ in $h$-balls. To convert this into an estimate for $Q$ in $g$-balls, we need the following distance comparison. Very similar results have previously been obtained in \cite{CatinoMastroliaRoncoroni} and \cite{ElbertNelliRosenberg}.

\begin{lemma}\label{distance comparison}
Assume $B = \overline{B_g(p,L)}$ is compactly contained in $\Sigma \setminus \partial \Sigma$, and that
\begin{equation}\label{eq:normalised}
Q \leq 4 \quad \text{on}\quad B.
\end{equation}
Set $D=\dist_h(p,\partial B)$, where $h=u^{4/3}g$. If $L \geq 4$, then we have
\[
L\leq C u(p)^{-8/3}D^4,
\]
where $C>0$ is a universal constant.
\end{lemma}
\begin{proof}
Let $\gamma : [0,T] \to \Sigma$ be an $h$-minimizing geodesic from $p$ to $\partial B$ which realises the distance $D$. Note that $\gamma$ remains in $B$; otherwise, truncating it at the first exit point yields a shorter geodesic to the boundary. We assume $g(\dot \gamma, \dot \gamma) = 1$, so that $T$ is the $g$-length of $\gamma$. It then follows that $T \geq L$. The $h$-length of $\gamma$ is given by 
\[D = \int_0^T \sqrt{h(\dot \gamma, \dot \gamma)} \,dt = \int_0^T e^{\psi(t)}\,dt, \qquad \psi(t) = \frac{2}{3}f(\gamma(t)).\]
We need to prevent $\psi$ from dropping too far below its initial value. 

Since $\gamma$ minimizes $h$-length, for any $\varphi : [0,T] \to \mathbb{R}$ of class $H_0^1$, we have the stability condition
\[
2\int_0^T e^{-\psi}\dot \varphi^2 \,dt \geq \int_0^T e^{-\psi}\Ric_h(\dot\gamma, \dot\gamma)\varphi^2 \,dt.
\]
This says that the $h$-Ricci curvature in the direction $\dot \gamma$ cannot be too large along $\gamma$. On the other hand, that same curvature component controls derivatives of $\psi$. Indeed, by \eqref{conformal Ricci} we have
\[
\Ric_h(\dot\gamma,\dot\gamma) = \Ric(\dot\gamma,\dot\gamma) - \nabla^2_h \psi(\dot\gamma,\dot\gamma) - \dot\psi^2 - \Delta \psi.
\]
Combining this with $\nabla^2_h \psi(\dot\gamma,\dot\gamma) + \dot\psi^2 = \ddot\psi$, which follows because $\gamma$ is a reparameterized $h$-geodesic, and the bounds 
\[
\Ric(\dot\gamma,\dot\gamma) = - |A(\dot\gamma,\cdot)|^2 \geq -\frac{2}{3}s, \qquad -\Delta \psi = \frac{2}{3}(q+s) \geq \frac{3}{2}\dot\psi^2 + \frac{2}{3}s,
\]
we obtain 
\[
\Ric_h(\dot\gamma,\dot\gamma) \geq -\ddot\psi + \frac{3}{2}\dot\psi^2.
\]
Therefore,
\[
2\int_0^T e^{-\psi}\dot \varphi^2 \,dt 
\geq \int_0^T e^{-\psi}\bigg(\frac{3}{2}\dot\psi^2 -\ddot\psi \bigg)\varphi^2 \,dt.
\]

It is convenient to clear the weight $e^{-\psi}$, by choosing a variation of the form $\varphi = e^{\psi/2}v$. After a short calculation, using the fundamental theorem of calculus to cancel
\[
2\dot\psi v\dot v + \ddot\psi v^2 = \frac{d}{dt}(\dot \psi v^2),
\]
we obtain the simpler inequality
\begin{equation}\label{stability}
2\int_0^T \dot v^2 \,dt \geq \int_0^T \dot\psi^2 v^2\,dt. 
\end{equation}

We use \eqref{stability} to bound the drop in $\psi$ by a constant times $\log L$. To achieve this, we choose $v$ so that
\[
v(t) = \sqrt{t}, \quad 1 \leq t \leq L/2,
\]
cut off linearly on the left and right,
\[
v(t) = t, \quad 0 \leq t \leq 1, \qquad v(t) = \sqrt{2/L}(L-t), \quad L/2 \leq t \leq L,
\]
and set $v=0$ on $[L,T]$. With this choice of $v$ we have 
\[
\int_0^T \dot v^2 \, dt = 2 + \frac{1}{4}\log(L/2)
\]
and hence, for $1 \leq t \leq L/2$, by \eqref{stability} and Cauchy--Schwarz,
\begin{align*}
\psi(t) - \psi(1) \geq -\bigg(\int_1^{L/2} t\dot\psi^2\,dt\bigg)^\frac{1}{2}
\bigg(\int_1^{L/2}\frac{1}{t}\,dt\bigg)^\frac{1}{2} \geq -6 - \frac{3}{4}\log(L/2).
\end{align*}
To bound $\psi(1)$ we recall that $\psi(0) = \frac{2}{3}\log u(p)$ and $Q \leq 4$ along $\gamma$. Given that 
\[
\dot \psi^2 = \frac{4}{9} \langle \nabla f, \dot \gamma\rangle^2 \leq \frac{4}{9}q \leq \frac{1}{4}Q \leq 1,
\]
the fundamental theorem of calculus yields $\psi(1) \geq \psi(0)-1$, and so we have
\[
\psi(t) \geq \psi(0) - 8 - \frac{3}{4}\log(L/2)
\]
for $1 \leq t \leq L/2$. 

Finally, we substitute the previous inequality into $D = \int_0^T e^{\psi(t)}\,dt$. This gives
\[
D \geq u(p)^{2/3}e^{-8}\bigg(\frac{L}{2} - 1\bigg)\bigg(\frac{L}{2}\bigg)^{-\frac{3}{4}}.
\]
Since $L\geq4$, the right-hand side is bounded from below by a universal constant times $u(p)^{2/3}L^{1/4}$. This completes the proof. 
\end{proof}

We now prove Theorem~\ref{main} by applying our bound for $R_h$ and the distance comparison lemma at a point where $\varphi Q$ is maximised, for a suitable cutoff $\varphi$. 

\begin{proof}[Proof of Theorem~\ref{main}]
Fix $p\in\Sigma\setminus\partial\Sigma$ and
$0<R<\dist_g(p,\partial\Sigma)$. Let $r$ denote the $g$-distance from $p$,
and choose $\bar p$ where $(R-r)^2Q$ attains its maximum
on $\overline{B_g(p,R)}$. Set
\[
\ell=\frac{R-r(\bar p)}{2},
\qquad
L=\ell\sqrt{Q(\bar p)}.
\]
Maximality gives
\[
R^2Q(p)\leq (R-r(\bar p))^2Q(\bar p)=4L^2,
\]
so it suffices to bound $L$ by a universal constant.
We may therefore assume $L\geq4$.

By the triangle inequality,
$B=\overline{B_g(\bar p,\ell)}$ is compactly contained in $B_g(p,R)$. Moreover,
we have $R-r(x)\geq\ell$ on $B$. Maximality of $\bar p$
therefore gives
\begin{equation*}\label{Q bound}
Q(x)\leq
\frac{(R-r(\bar p))^2}{(R-r(x))^2}Q(\bar p)
\leq4Q(\bar p)
\end{equation*}
for each $x \in B$.

Dilate the immersion by $\sqrt{Q(\bar p)}$. In the rescaled geometry,
$B$ has $g$-radius $L$, and
\[
Q(\bar p)=1,
\qquad
Q\leq4\quad\text{on }B.
\]
We set
\[
D=\dist_h(\bar p,\partial B),
\qquad h=u^{4/3}g.
\]
Theorem~\ref{conf max princ} then gives
\[
u(\bar p)^{-4/3}D^2=R_h(\bar p)D^2\leq K,
\]
and hence by Lemma~\ref{distance comparison} we have
\[
L\leq C u(\bar p)^{-8/3}D^4
=C\big(R_h(\bar p)D^2\big)^2\leq CK^2.
\]

Returning to the original scale, we conclude that
\[
R^2Q(p)\leq4L^2\leq 4C^2 K^4.
\]
Letting $R\uparrow\dist_g(p,\partial\Sigma)$ and using
$q+s\leq\frac35Q$ completes the proof.
\end{proof}

\iffalse
Let us comment on our choice $k = \frac{4}{3}$. In fact all of our arguments also go through for $\frac{4}{3} < k < 2$. The main difference for these other values is that the Bakry--\'{E}mery~Ricci curvature of $h = u^k g$ is no longer nonnegative, but only admits a lower bound in terms of $W$. This turns out to be sufficient---an extra curvature term appears in the localisation step, but this term is of lower order and is easily absorbed. 

There is also some freedom in how we choose the quantity $Q$. For $Q = q + c s$ with $c > 0$, the calculation which established Proposition~\ref{Delta Q} instead yields
\begin{equation}\label{Delta Q general}
\Delta Q\geq
\frac34\frac{|\nabla Q|^2}{Q}
-\langle\nabla f,\nabla Q\rangle
+q^2
+\left(\frac{20}{3}-3c-\frac{10}{3c}\right)qs
+\left(\frac23-2c\right)s^2
\end{equation}
at every point where $Q>0$. Our choice $c=\frac{2}{3}$ is convenient
because it reduces the reaction terms to $(q-s)Q$. We have recorded the general inequality \eqref{Delta Q general} here because, in other contexts, a different choice of $c$ may be preferred. 
\fi

\section{An alternative proof in dimension two}\label{sec 2D proof}

The previous two sections established Theorem~\ref{main} for $n = 3$. This implies the $n=2$ case---take the product with a line. Here we record an alternative proof in dimension two, which does not require us to make a conformal change of metric. 

Let $\Sigma \to \mathbb{R}^3$ be a minimal immersion, and let $u:\Sigma\to\mathbb{R}$ be a positive Jacobi field on $\Sigma$. We continue to use the notation
\[
f = \log u, \qquad q = |\nabla f|^2, \qquad s = |A|^2.
\]

\begin{proposition}\label{2D Delta Q}
Let $c$ be a positive constant and set $Q=q+cs$. At every point where $Q > 0$ we have
\begin{equation}\label{eq: 2D Delta Q}
\Delta Q \geq \frac{|\nabla Q|^2}{Q} + 2q^2 + \bigg(3-\frac{1}{c}\bigg)qs + (1-2c)s^2.
\end{equation}
\end{proposition}
\begin{proof}
We first consider a point where both $q$ and $s$ are positive. Combining the Bochner formula with $\Delta f = -s-q$ and $\Ric = -A^2 = -\frac{1}{2} s g$, one obtains
\[
\Delta q = 2|\nabla^2 f|^2 -2\langle\nabla f, \nabla s + \nabla q\rangle - sq.
\]
To rewrite the Hessian term, we choose an orthonormal frame with $e_1 = \nabla f/|\nabla f|$ and use
\[
\nabla_i q = 2\nabla^2 f(\nabla f, e_i) = 2|\nabla f|f_{1i}, \qquad f_{22} = \Delta f - f_{11}.
\]
This leads to
\begin{equation}\label{2D Hessian}
2|\nabla^2 f|^2 - 2\langle\nabla f, \nabla q\rangle = \frac{|\nabla q|^2}{q} + 2s\langle\nabla f, \nabla \log q\rangle  + 2(q+s)^2,
\end{equation}
and hence
\[
\Delta q = \frac{|\nabla q|^2}{q} + 2s\langle\nabla f, Y\rangle + 2q^2 + 3qs +2s^2,
\]
where we have set $Y = \nabla \log q - \nabla \log s$.

In dimension two, Simons' identity is
\[
\Delta s = \frac{|\nabla s|^2}{s} - 2s^2,
\]
and so for $Q = q + cs$ we have
\[
\Delta Q = \frac{|\nabla q|^2}{q} + c\frac{|\nabla s|^2}{s} + 2s\langle\nabla f, Y\rangle + 2q^2 + 3qs + 2(1-c)s^2.
\]

A short calculation shows that
\[
\frac{|\nabla q|^2}{q} + c\frac{|\nabla s|^2}{s} = \frac{|\nabla Q|^2}{Q} + c\frac{sq}{Q}|Y|^2.
\]
Moreover, completing the square gives
\[
c\frac{sq}{Q}|Y|^2 + 2s\langle\nabla f, Y\rangle = c\frac{sq}{Q}\bigg|Y + \frac{Q}{cq}\nabla f\bigg|^2 - \frac{sQ}{c}.
\]
Combining these facts, and throwing away a nonnegative term, we arrive at the required estimate:
\[
\Delta Q \geq \frac{|\nabla Q|^2}{Q} + 2q^2 + \bigg(3-\frac{1}{c}\bigg)qs + (1-2c)s^2.
\]

It remains to consider points where $Q > 0$, but either $q$ or $s$ vanishes. At a point where $q = 0$, we have $s > 0$ and $\nabla q = 0$. Combining the Bochner formula, Simons' identity and the trace inequality
\[
2|\nabla^2 f|^2 \geq (\Delta f)^2 = s^2,
\]
we obtain
\[
\Delta Q = 2|\nabla^2 f|^2 + c\frac{|\nabla s|^2}{s} - 2cs^2 \geq c\frac{|\nabla s|^2}{s} + (1-2c)s^2 = \frac{|\nabla Q|^2}{Q} + (1-2c)s^2.
\]
If instead $s = 0$, then $q > 0$ and $\nabla s = 0$. In this case Simons' identity implies $\Delta s \geq 0$, and combining Bochner with \eqref{2D Hessian} gives
\[
\Delta Q \geq \frac{|\nabla q|^2}{q} + 2q^2 \geq \frac{|\nabla Q|^2}{Q} + 2 q^2.
\]
In both cases we obtain the required inequality. 
\end{proof}

We now describe how Proposition~\ref{2D Delta Q} can be used to obtain Theorem~\ref{main} in dimension two. Setting $c = \frac{1}{3}$, so that $Q = q + \frac{1}{3} s$, \eqref{eq: 2D Delta Q} implies
\[
\Delta Q \geq \frac{|\nabla Q|^2}{Q} +  \frac{6}{5}Q^2.
\]
Suppose $\Sigma$ is complete and fix $p\in\Sigma\setminus\partial\Sigma$ and $0<R<\dist_g(p,\partial\Sigma)$. We introduce a cutoff of the form $\varphi(x) = (R^2 - r(x)^2)^2$, where $r(x)$ is the $g$-distance from $p$. Writing $\zeta=R^2-r^2$, at a positive interior maximum $\bar p$ for the function $\varphi Q$, one finds that 
\[
-4\zeta(1+r\Delta r)-8r^2+\frac{6}{5}\varphi Q\leq0.
\]
At cut points, this is made rigorous by smooth approximation of $r$ through upper barriers (see Appendix~\ref{Laplace comparison}).  Along a minimizing geodesic from $p$ to $\bar p$, we have $\varphi\geq\varphi(\bar p)$. Maximality of $\varphi Q$ therefore gives $Q\leq Q(\bar p)$ along this geodesic. Since $s\leq3Q$, the Ricci tensor satisfies
\[
\Ric = - \frac{1}{2}s g \geq -\frac{3}{2}Q(\bar p) g,
\]
and hence the classical Laplacian comparison theorem\footnote{The exact version we need is Proposition~\ref{weighted Laplacian comparison} with $N = n$ and $\psi$ constant.} implies
\[
\Delta r(\bar p) \leq \sqrt{\frac{3}{2}Q(\bar p)} \coth\bigg(\sqrt{\frac{3}{2} Q(\bar p)}r(\bar p)\bigg) \leq C\bigg(\sqrt{Q(\bar p)} + \frac{1}{r(\bar p)}\bigg)
\]
for some universal constant $C$. Using this estimate, we conclude that 
\[
-CR^2 - CR\sqrt{\varphi(\bar p) Q(\bar p)} +  \frac{6}{5}\varphi(\bar p) Q(\bar p) \leq 0.
\]
The upshot is that Laplacian comparison leads only to lower-order terms. Using Young's inequality and rearranging gives $R^2Q(p)\leq C$. Letting $R\uparrow\dist_g(p,\partial\Sigma)$ and using $q+s\leq3Q$ completes the proof.

\appendix

\section{An identity for the conformal scalar curvature}\label{sharp BS}

We return to the three-dimensional setting and notation of Section~\ref{sec BS}. With $h=u^{4/3}g$, we split the scalar curvature $R_h$ into the two contributions
\[
X=\frac{16}{9}u^{-4/3}q,\qquad
Y=\frac53u^{-4/3}s,\qquad R_h=X+Y.
\]
In this appendix we refine Proposition~\ref{Jacobi scalar inequality}, keeping track of the nonnegative terms which were previously discarded. 

The Codazzi equations and minimality imply that $\nabla A$ lies in the space of $(0,3)$-tensors which are totally symmetric and trace-free. We denote by $B$ the orthogonal projection of $df \otimes A$ onto that same space. 

\begin{proposition}\label{scalar defect identity}
The scalar curvature of $h$ satisfies
\begin{equation}\label{full scalar defect}
\Delta_h R_h
=\frac18X^2+\frac{767}{1000}XY+\frac{2}{75}Y^2+\mathcal D,
\end{equation}
where
\begin{align}\label{scalar defect squares}
u^{8/3}\mathcal D= \frac{32}{9}|T|^2
+\frac{10}{3}\left|\nabla A-\frac{31}{15}B\right|^2 +\frac{6244}{675}\left(\frac23sq-|A(\nabla f,\cdot)|^2\right).
\end{align}
The function $\mathcal D$ is nonnegative.
\end{proposition}
\begin{proof}
The definition of $B$ gives the first identity below, and a direct calculation gives the second:
\[
\langle\nabla A,B\rangle=\frac12\langle\nabla f,\nabla s\rangle,
\qquad
|B|^2=\frac13qs+\frac25|A(\nabla f,\cdot)|^2.
\]
Completing the square, we obtain
\begin{align*}
\frac{10}{3}|\nabla A|^2
-\frac{62}{9}\langle\nabla f,\nabla s\rangle
&=\frac{10}{3}\left|\nabla A-\frac{31}{15}B\right|^2
-\frac{1922}{135}|B|^2.
\end{align*}
Substituting this into \eqref{Bochner Simons scalar} and using the formula for $|B|^2$ gives
\[
u^{8/3}\big(\Delta_h R_h-\mathcal D\big)
=\frac{32}{81}q^2+\frac{1534}{675}qs+\frac{2}{27}s^2.
\]
The definitions of $X$ and $Y$ now give \eqref{full scalar defect}. Since $A$ is traceless, we have $|A(\nabla f,\cdot)|^2\leq\frac23qs$, so every term in $\mathcal D$ is nonnegative.
\end{proof}

\section{Weighted Laplacian comparison}\label{Laplace comparison}

We record a local form of Qian's weighted Laplacian comparison theorem~\cite{Qian}. Qian works with a complete metric and a global lower bound for its finite-dimensional Bakry--\'{E}mery Ricci tensor. The proof only uses the curvature bound along the minimizing geodesic, so the following local form does not require completeness.

\begin{proposition}\label{weighted Laplacian comparison}
Let $(M,h)$ be a Riemannian manifold of dimension $n\geq2$, possibly with boundary, and let $\psi$ be a smooth function on its interior. Let $\gamma$ be a unit-speed minimizing geodesic from $p$ to $x$, which lies in the interior of $M$ and has length $\ell>0$. If, for some $N>n$ and $\kappa\geq0$,
\[
\bigg(\Ric_h+\nabla^2_h\psi
-\frac{1}{N-n}d\psi\otimes d\psi\bigg)
(\dot\gamma,\dot\gamma)\geq -(N-1)\kappa^2
\]
holds along $\gamma$, then the distance function $r=\dist_h(p,\cdot)$ satisfies
\[
\Delta_{h,\psi} r(x)
\leq (N-1)\kappa\coth\bigl(\kappa r(x)\bigr)
\]
in the barrier sense, where $\Delta_{h,\psi}
=\Delta_h-\langle\nabla^h\psi,\nabla^h\cdot\rangle_h$.
When $\kappa=0$, the right-hand side is interpreted as $(N-1)/r(x)$.
If $\psi$ is constant, we also allow $N=n$, with the last term in the curvature hypothesis omitted.
\end{proposition}

\begin{proof}
Suppose first that $N>n$ and $r$ is smooth near $x$. For any smooth function $v:[0,\ell]\to\mathbb{R}$ with $v(0)=0$ and $v(\ell)=1$, the second variation of length gives
\[
\Delta_h r(x)
\leq\int_0^\ell
\big((n-1)\dot v^{\,2}
-v^2\Ric_h(\dot\gamma,\dot\gamma)\big)\,dt.
\]
Using the curvature hypothesis and integrating the term $v^2\ddot\psi$ by parts, we obtain
\begin{align*}
\Delta_{h,\psi} r(x)
&\leq\int_0^\ell
\bigg((n-1)\dot v^{\,2}
-2v\dot v\dot\psi
-\frac{v^2\dot\psi^{\,2}}{N-n}
+(N-1)\kappa^2v^2\bigg)\,dt\\
&=(N-1)\int_0^\ell
\big(\dot v^{\,2}+\kappa^2v^2\big)\,dt
-\frac{1}{N-n}\int_0^\ell
\big(v\dot\psi+(N-n)\dot v\big)^2\,dt.
\end{align*}
Taking
\[
v(t)=\frac{\sinh(\kappa t)}{\sinh(\kappa\ell)}
\]
when $\kappa$ is nonzero, and $v(t)=t/\ell$ when $\kappa$ is zero, gives
\[
\int_0^\ell\big(\dot v^{\,2}+\kappa^2v^2\big)\,dt
=[v\dot v]_0^\ell
=\kappa\coth(\kappa\ell).
\]
This proves the claim. 

To remove the smoothness assumption, we use Calabi's trick \cite{Calabi}. Fix $0<\varepsilon<\ell$. Since $\gamma$ is minimizing, $x$ is not conjugate to $\gamma(\varepsilon)$ along the shortened segment. Varying its endpoint therefore gives a smooth family of geodesics from $\gamma(\varepsilon)$ to points $y$ near $x$. Let $r_\varepsilon(y)$ be $\varepsilon$ plus the length of the corresponding geodesic. Then $r_\varepsilon$ is smooth near $x$, and
\[
r_\varepsilon\geq r,
\qquad
r_\varepsilon(x)=r(x).
\]
The same stability calculation, applied along the shortened segment, gives
\[
\Delta_{h,\psi} r_\varepsilon(x)
\leq (N-1)\kappa\coth\bigl(\kappa(\ell-\varepsilon)\bigr).
\]
Letting $\varepsilon\downarrow0$ proves the upper-barrier inequality. Since each shortened segment lies on $\gamma$, the argument requires no curvature bound away from $\gamma$.
The case $N=n$ with constant $\psi$ follows by letting $N\downarrow n$.
\end{proof}

\section{Conformal transformation formulae}
\label{conformal formulae}

Let $g$ be a Riemannian metric on an $n$-dimensional manifold and set $h=e^{2\phi}g$. We record here the conformal transformation formulae made use of in this article. All unadorned derivatives, norms and inner products refer to $g$. For any smooth
function $v$, we have
\begin{align*}
\nabla^h v
&=e^{-2\phi}\nabla v,\\
%|\nabla^h v|_h^2
%&=e^{-2\phi}|\nabla v|^2,\\
\nabla_h^2v
&=\nabla^2v-d\phi\otimes dv-dv\otimes d\phi
+\langle\nabla\phi,\nabla v\rangle g,\\
\Delta_hv
&=e^{-2\phi}\big(
\Delta v+(n-2)\langle\nabla\phi,\nabla v\rangle
\big),\\
\Ric_h
&=\Ric_g-(n-2)\big(\nabla^2\phi-d\phi\otimes d\phi\big)
-\big(\Delta\phi+(n-2)|\nabla\phi|^2\big)g\\
%&=\Ric_g-(n-2)\nabla_h^2\phi
%-(n-2)d\phi\otimes d\phi-\Delta\phi\,g,\\
R_h
&=e^{-2\phi}\big(
R_g-2(n-1)\Delta\phi
-(n-1)(n-2)|\nabla\phi|^2
\big).
\end{align*}


\begin{thebibliography}{99}

\bibitem{AntonelliXu}
G.~Antonelli and K.~Xu,
\emph{New spectral Bishop--Gromov and Bonnet--Myers theorems and applications to isoperimetry},
to appear in J. Eur. Math. Soc.,
\href{https://arxiv.org/abs/2405.08918}{arXiv:2405.08918}.

\bibitem{Bellettini}
C.~Bellettini,
\emph{Extensions of Schoen--Simon--Yau and Schoen--Simon theorems via iteration \`a la De Giorgi},
Invent. Math. \textbf{240} (2025), 1--34.

\bibitem{Bernstein}
S.~Bernstein,
\emph{Sur un th\'eor\`eme de g\'eom\'etrie et son application aux \'equations aux d\'eriv\'ees partielles du type elliptique},
Comm. Soc. Math. Kharkow (2) \textbf{15} (1915), 38--45.

\bibitem{BombieriDeGiorgiGiusti}
E.~Bombieri, E.~De Giorgi and E.~Giusti,
\emph{Minimal cones and the Bernstein problem},
Invent. Math. \textbf{7} (1969), 243--268.

\bibitem{BombieriGiusti}
E.~Bombieri and E.~Giusti,
\emph{Harnack's inequality for elliptic differential equations on minimal surfaces},
Invent. Math. \textbf{15} (1972), 24--46.

\bibitem{CabreCatinoMariMastroliaRoncoroni}
X.~Cabr\'e, G.~Catino, L.~Mari, P.~Mastrolia and A.~Roncoroni,
\emph{Gradient estimates for the Green kernel under spectral Ricci bounds, and the stable Bernstein theorem in $\mathbb{R}^4$},
preprint (2026), \href{https://arxiv.org/abs/2604.14393}{arXiv:2604.14393}.

\bibitem{CaffarelliNirenbergSpruck}
L.~Caffarelli, L.~Nirenberg and J.~Spruck,
\emph{On a form of Bernstein's theorem},
in \emph{Analyse math\'ematique et applications},
Gauthier-Villars, Montrouge, 1988, pp.~55--66.

\bibitem{Calabi}
E.~Calabi,
\emph{An extension of E.~Hopf's maximum principle with an
application to Riemannian geometry},
Duke Math. J. \textbf{25} (1958), 45--56.

\bibitem{CatinoMastroliaRoncoroni}
G.~Catino, P.~Mastrolia and A.~Roncoroni,
\emph{Two rigidity results for stable minimal hypersurfaces},
Geom. Funct. Anal. \textbf{34} (2024), 1--18.

\bibitem{ChengYau}
S.-Y.~Cheng and S.-T.~Yau,
\emph{Differential equations on Riemannian manifolds and their geometric applications},
Comm. Pure Appl. Math. \textbf{28} (1975), 333--354.

\bibitem{ChodoshLiAnisotropic}
O.~Chodosh and C.~Li,
\emph{Stable anisotropic minimal hypersurfaces in $\mathbb{R}^4$},
Forum Math. Pi \textbf{11} (2023), Paper No.~e3, 22~pp.

\bibitem{ChodoshLi}
O.~Chodosh and C.~Li,
\emph{Stable minimal hypersurfaces in $\mathbb{R}^4$},
Acta Math. \textbf{233} (2024), 1--31.

\bibitem{ChodoshLiMinterStryker}
O.~Chodosh, C.~Li, P.~Minter and D.~Stryker,
\emph{Stable minimal hypersurfaces in $\mathbb{R}^5$},
Ann. of Math. (2) \textbf{204} (2026), 423--453.

\bibitem{Colding}
T.~H.~Colding,
\emph{New monotonicity formulas for Ricci curvature and applications. I},
Acta Math. \textbf{209} (2012), 229--263.

\bibitem{doCarmoPeng}
M.~do Carmo and C.~K.~Peng,
\emph{Stable complete minimal surfaces in $\mathbb{R}^3$ are planes},
Bull. Amer. Math. Soc. (N.S.) \textbf{1} (1979), 903--906.

\bibitem{EckerHuisken}
K.~Ecker and G.~Huisken,
\emph{A Bernstein result for minimal graphs of controlled growth},
J. Differential Geom. \textbf{31} (1990), 397--400.

\bibitem{ElbertNelliRosenberg}
M.~F.~Elbert, B.~Nelli and H.~Rosenberg,
\emph{Stable constant mean curvature hypersurfaces},
Proc. Amer. Math. Soc. \textbf{135} (2007), 3359--3366.

\bibitem{FischerColbrie}
D.~Fischer-Colbrie,
\emph{On complete minimal surfaces with finite Morse index
in three manifolds},
Invent. Math. \textbf{82} (1985), 121--132.

\bibitem{FischerColbrieSchoen}
D.~Fischer-Colbrie and R.~Schoen,
\emph{The structure of complete stable minimal surfaces in $3$-manifolds of nonnegative scalar curvature},
Comm. Pure Appl. Math. \textbf{33} (1980), 199--211.

%\bibitem{GulliverLawson}
%R.~Gulliver and H.~B.~Lawson, Jr.,
%\emph{The structure of stable minimal hypersurfaces near a singularity},
%in \emph{Geometric Measure Theory and the Calculus of Variations},
%Proc. Sympos. Pure Math., vol.~44,
%Amer. Math. Soc., Providence, RI, 1986, pp.~213--237.

\bibitem{Hamilton}
R.~S.~Hamilton,
\emph{Harnack estimate for the mean curvature flow},
J. Differential Geom. \textbf{41} (1995), 215--226.

\bibitem{HardtSimon}
R.~Hardt and L.~Simon,
\emph{Area minimizing hypersurfaces with isolated singularities},
J. Reine Angew. Math. \textbf{362} (1985), 102--129.

\bibitem{HongLiWang}
H.~Hong, H.~Li and G.~Wang,
\emph{The stable Bernstein theorem in $\mathbb{R}^7$},
preprint (2026), \href{https://arxiv.org/abs/2609.15720}{arXiv:2609.15720}.

\bibitem{LiYau}
P.~Li and S.-T.~Yau,
\emph{On the parabolic kernel of the Schr\"odinger operator},
Acta Math. \textbf{156} (1986), 153--201.

\bibitem{Mazet}
L.~Mazet,
\emph{Stable minimal hypersurfaces in $\mathbb{R}^6$},
preprint (2024), \href{https://arxiv.org/abs/2405.14676}{arXiv:2405.14676}.

\bibitem{Pogorelov}
A.~V.~Pogorelov,
\emph{On the stability of minimal surfaces},
Soviet Math. Dokl. \textbf{24} (1981), 274--276.

\bibitem{Qian}
Z.~Qian,
\emph{Estimates for weighted volumes and applications},
Quart. J. Math. Oxford Ser. (2) \textbf{48} (1997), 235--242.

\bibitem{Schoen}
R.~Schoen,
\emph{Estimates for stable minimal surfaces in three-dimensional manifolds},
in \emph{Seminar on Minimal Submanifolds},
Ann. of Math. Stud., vol.~103, Princeton Univ. Press, Princeton, NJ, 1983, pp.~111--126.

\bibitem{SchoenSimon}
R.~Schoen and L.~Simon,
\emph{Regularity of stable minimal hypersurfaces},
Comm. Pure Appl. Math. \textbf{34} (1981), 741--797.

\bibitem{SchoenSimonYau}
R.~Schoen, L.~Simon and S.-T.~Yau,
\emph{Curvature estimates for minimal hypersurfaces},
Acta Math. \textbf{134} (1975), 275--288.

\bibitem{SchoenYauStructure}
R.~Schoen and S.-T.~Yau,
\emph{On the structure of manifolds with positive scalar curvature},
Manuscripta Math. \textbf{28} (1979), 159--183.

\bibitem{Simons}
J.~Simons,
\emph{Minimal varieties in Riemannian manifolds},
Ann. of Math. (2) \textbf{88} (1968), 62--105.

\bibitem{Yau}
S.-T.~Yau,
\emph{Harmonic functions on complete Riemannian manifolds},
Comm. Pure Appl. Math. \textbf{28} (1975), 201--228.

\end{thebibliography}
\end{document}